\documentclass[11pt]{article}
\usepackage{layout}
\usepackage{amssymb}
\usepackage{epsfig}
\usepackage{graphicx}
\usepackage{url}
\usepackage{color}
\usepackage{latexsym}
\usepackage[latin1]{inputenc}
\usepackage{verbatim}

\usepackage[left=2cm,right=2cm,top=2cm,bottom=2cm]{geometry}
\usepackage{tikz}

\usepackage{amssymb}
\usepackage{algorithmicx}
\usepackage[ruled]{algorithm}
\usepackage{algpseudocode}
\usepackage{algc}
\usepackage{multirow}
\alglanguage{pseudocode}
\usepackage{graphicx}
\usepackage{amsmath,epsfig,graphics}
\usepackage{arydshln}  

\newtheorem{theorem}{Theorem}[section]

\newtheorem{remark}[theorem]{Remark}

\newcommand{\R}{\mathbb{R}}

\title{Derivative-free optimization approach for structured symmetric matrices with fixed eigenvalues}

\author{Carmo P. Br\'as\thanks{Center for Mathematics and Applications (NOVA Math) and Department of Mathematics,
		NOVA School of Science and Technology (NOVA FCT), 2829-516 Caparica, Portugal ({\tt mb@fct.unl.pt}).}
	\and Evelin H. M. Krulikovski\thanks{Department of Mathematics, Federal University of Paran\'a, Curitiba,
		Brazil ({\tt evelin.krulikovski@ufpr.br}).}
	\and Marcos Raydan\thanks{Center for Mathematics and Applications (NOVA Math) and Department of Mathematics,
		NOVA School of Science and Technology (NOVA FCT), 2829-516 Caparica, Portugal  ({\tt m.raydan@fct.unl.pt}).}}

\begin{document}
	
\date{May 19, 2026}

\maketitle

\begin{abstract}
	A Derivative-Free Optimization (DFO) model is developed and analyzed for solving inverse structured symmetric
	matrix problems for which the eigenvalues are specified. Some (zero and nonzero) entries are preassigned and
	cannot be changed, while others should be nonzero but their values are not given.
	The rest of the entries are completely free. The obtained matrix must meet these criteria and have the
	 specified eigenvalues. This specialized inverse eigenvalue problem is relevant to various applications
	 and is linked to determining the graph, with weights on the undirected edges, of the matrix associated
	 with its sparse pattern. Our novel optimization model requires computing the eigenvalues of a
	 symmetric matrix to evaluate the non-differentiable objective function.
	 We apply   deterministic DFO schemes, specifically
	the global variant GLODS of the well-known family of directional direct search (DDS) methods.
  We discuss its convergence properties  which are based on the fact that the objective function of our model
  is Lipschitz continuous.  Additionally, we explore the potential benefits of using several well-established
   heuristic strategies  to solve the proposed optimization model. We present preliminary
   	 numerical results to illustrate   and compare the performance of the considered deterministic
   	  and  heuristic DFO options in various possible scenarios.
 \end{abstract}
{\small
	{\bf Keywords:} {Inverse eigenvalue problems, spectral graph theory, derivative-free optimization,
	 global optimization, bound constraints, heuristic strategies  }

\section{Introduction} \label{Intro}

 Given a set of real eigenvalues in increasing order, say $\{\lambda_1, \lambda_2,\dots,\lambda_n\}$,
 that defines a fixed diagonal
  matrix $\Lambda$ whose diagonal entries are the given $\lambda_i$s, the problem is to find a real symmetric
   $n\times n$ matrix $A$ such that its eigenvalues are the given ones. \\

    There are some additional and fundamental requirements that add difficulty to the problem. For some specified
    pairs $(i,j)\in I_{0}$ the entries  $a_{ij}$ of $A$ must be zero at the solution. For some other specified
     pairs $(i,j)\in I_{fix}$ the entries  $a_{ij}$ of $A$ must have some given nonzero values  at the solution.
  For some other given pairs $(i,j)\in I_{nz}$ the entries  $a_{ij}$ of $A$ are free  but they must
  	 be bounded away from zero by some given thresholds.  Finally, for some other specified
    pairs $(i,j)\in I_{free}$ the entries  $a_{ij}$ of $A$ are free variables that can take any real value.
    The  four given index sets: $I_{0}, I_{fix}, I_{nz}$ and $I_{free}$ should conform a partition
     of all the possible pairs $(i,j)$, with $i$ and $j$ in  $\{1,2,\dots, n\}$ and they need to
     reflect the symmetry of the desired matrix $A$. Some examples and  some additional details
     describing this problem, as well as a  manifold optimization approach (based on off-the-shelf packages)
    for low-dimensional problems can be found in \cite{sutton}. More recently, a penalization strategy, an
     augmented Lagrangian approach, and a block alternating technique, all of which are  based on manifold
      optimization, have been introduced and analyzed in \cite{chehab}.  It should be noted that all of the
       aforementioned methods use the $n \times n$ orthogonal matrix of the eigenvectors of the desired solution
    as the only variable.   The model proposed and solved in \cite{chehab} includes a key term
  	in the objective function that ensures that the sum of all elements in $I_{nz}$ is as far from zero as possible.
     However, since this term only depends on the implicit orthogonal matrix of eigenvectors and does not
     treat the matrix's entries as variables, it cannot guarantee that every entry in $I_{nz}$  is sufficiently
     far from zero.  \\

    The desired  pattern, described above for the unknown symmetric matrix $A$, can be associated with a weighted
     undirected graph with $n$ vertices,     in which two vertices $i$ and $j$ are connected (adjacent) if and
      only if $a_{ij}$ is not zero. The described inverse problem and the connection between  the undirected graph
      and the matrix $A$ with preset eigenvalues appear in several real-life applications related to quantum
      chemistry,  information theory, social networks, and combinatorial optimization; see, e.g.,
      \cite{Cvet09, Chu, Trinaj, Wasser} and the references cited therein.
    The bounds (or thresholds) imposed on the weights of adjacent vertices associated with the
   	 $I_{nz}$ set appear naturally in these applications. They ensure that the weights of those edges in the
   	 solution graph are sufficiently far from zero.
    It is worth noting that not every arbitrary assignment of real eigenvalues can guarantee the existence of a
     certain structure of the associated symmetric matrix. This is a topic of recent interest in matrix analysis.
      For a unified theoretical development of whether a predetermined distribution of eigenvalues can be
      achieved with a given graph, and also how the graph of a symmetric matrix influences the possible
       multiplicities of its  eigenvalues, we refer to the book by Johnson and Saiago \cite{JohnSaiago}. \\

  To address these problems of inverse structured matrices, we propose a new optimization model with simple
   bound constraints whose non-differentiable objective function requires the computation of the eigenvalues
    of a symmetric matrix.    The simple bounds are imposed on the entries in the
    $I_{nz}$ index set to guarantee that the given thresholds are satisfied. We also describe in detail some
 iterative schemes already existing in the DFO literature (deterministic and heuristic) to solve it;
  see, e.g., \cite{AudetHare, Chopard, ConnEtAl}.  \\

 The remainder of this paper is organized as follows.  {In Section \ref{dfoapp}, we introduce the
 notation and  key  components required to accurately present the proposed optimization model.  Furthermore,
  we examine in detail  the model's objective function, with particular emphasis on its smoothness properties.}
     In Section \ref{dds},   we  briefly describe the Directional Direct Search (DDS) family of DFO
     methods \cite{AudetHare, ConnEtAl},   focusing on the global variant: Global and Local
      Optimization using Direct Search (GLODS) \cite{Custodio15}. We also discuss
     the convergence properties of  GLODS {and how the existing theory applies to the proposed
     	 optimization model.}   In Section \ref{Heurist}, we briefly describe some well-established heuristic
      strategies that can also be applied to our optimization problem. {In Section \ref{NumRes}, we report the numerical results obtained on various test problems
   	and provide additional insights into the deterministic and heuristic schemes considered. We also compare these
   DFO   schemes with the recent differentiable manifold optimization approach \cite{chehab}. Numerous numerical
       experiments on matrices of various sizes and patterns demonstrate the effectiveness of the proposed DFO
    approach.}  Finally, in Section \ref{conclu}, we present some final comments and perspectives.

\section{Derivative-free optimization approach} \label{dfoapp}

 For any arbitrary symmetric matrix, say $X_0$, that satisfies the requirements associated with the
 index sets $I_{0}$  and $I_{fix}$, and that contains (say randomly) some real values for the entries
 in the index set  $I_{free}$, and some nonzero values (random but within the given bounds) for the
  entries in the index  set $I_{nz}$, we can perform the following computational steps.
   Compute the eigenvalues of $X_0$,  using any available software package, and obtain a diagonal
   matrix   $\Lambda(X_0)$ whose diagonal entries are the eigenvalues of $X_0$ in increasing order.
  Then, compare the obtained eigenvalues  with the given matrix $\Lambda$, by evaluating
  $\|vec(\Lambda(X_0) - \Lambda)\|_2$, where
 the operation $vec$ of a diagonal matrix consists in building an $n$-dimensional vector whose entries are the
  diagonal  entries of the argument matrix. If the 2-norm of that $n$-dimensional vector is zero, we can
  conclude that we were extremely fortunate, and that $X_0$ is a solution to our inverse matrix problem.
  Otherwise,  we have to try again through an iterative process. \\

   Roughly speaking, our intent is to produce a sequence   of matrices $X_k$, starting from $X_0$,
    that satisfy the specified
  structural requirements  and that drives  $\|vec(\Lambda(X_k) - \Lambda)\|_2$ down to zero, where
 $\Lambda(X_k)$ is a diagonal matrix containing the eigenvalues of $X_k$ in increasing order.  To be more specific,
 the main objective of our approach is to apply a DFO iterative scheme to produce a sequence of
 matrices that satisfy the
  structural requirements while minimizing a suitable objective function which includes the term
   $\|vec(\Lambda(X_k) - \Lambda)\|_2$.
 We note that  the use of some convenient DFO scheme is essential since the partial derivatives of
  $\|vec(\Lambda(X_k) - \Lambda)\|_2$ with respect to the  entry-wise variables $x_{ij}$
  (such that $(i,j) \in I_{nz}\cup I_{free}$)  exist only in very special cases,
  and even then they are not available. \\

   Concerning the objective function, it must be enriched to penalize the
    entries in the  set $ I_{nz}$ to keep them as far as
    possible from zero. In that sense the objective function to be minimized using DFO strategies is:
    \begin{equation} \label{ObjF}
      f(X) = \tau \;\|vec(\Lambda(X) - \Lambda)\|_2 - \sum_{(i,j)\in I_{nz}}  \ln(|x_{ij}|)\;,
       \end{equation}
  where $\tau>0$  is a predetermined fixed parameter that depends on the characteristics of the problem,
  and $\Lambda(X)$ is obtained for any possible candidate $X$, as described above.
 Note that the variables of this optimization problem are only the entries of the matrix $X$
  associated with the indices $I_{nz}$ and $I_{free}$. Therefore, the simple-bound constrained optimization
   problem to be considered  is
    \begin{equation} \label{optP}
   	\min_{x_{ij}:\:(i,j) \in I_{nz}\cup I_{free}} \hspace{-2mm} f(X) \quad \textrm{ subject to } \quad
   	l_{ij} \leq	x_{ij} \leq u_{ij}, 
   \end{equation}
  where  $l_{ij}$ and  $u_{ij}$ (for all  $(i,j) \in I_{nz}$) are provided  in advance and depend on
  the specific  inverse matrix problem being solved. If for a certain entry $(i,j) \in I_{nz}$ it is
  required that the final value of $a_{ij}\leq -\rho$, for some given $\rho>0$, then we set
  $l_{ij}=-M$ and $u_{ij}= -\rho$, where $0<M<+\infty$ is a sufficiently large positive real number.
  On the other hand, if for a certain entry $(i,j) \in I_{nz}$ it is
  required that the final value of $a_{ij}\geq \rho>0$,  then we set  $l_{ij}=\rho$ and $u_{ij}= +M$.
    For the free variables, that is the entries $(i,j) \in I_{free}$,   we set $l_{ij}=-M$
  	and $u_{ij}= +M$.
   Note that due to the nature of the
    bounds imposed on the $x_{ij}$ variables   (for all  $(i,j) \in I_{nz}\cup I_{free}$), the
    function $f(X)$ is bounded-below on the feasible region $\Omega$, which is formally defined as
   \begin{equation} \label{Omega}
     \Omega = \{X\in \R^{n\times n}\: :\:  -\infty < l_{ij} \leq	x_{ij} \leq u_{ij}< +\infty \;\; \textrm{ for }\; (i,j) \in I_{nz}\cup I_{free}\}.
   \end{equation}
   We also note that the penalization parameter $\tau>0$ in (\ref{ObjF}) plays a key
role to drive the term $\|vec(\Lambda(X) - \Lambda)\|_2$ down to zero, and
should be chosen depending on the specific characteristics of the given inverse matrix problem.
 {It must certainly be large enough to prioritize the minimization of the eigenvalue mismatch
 	  term ($\|vec(\Lambda(X) - \Lambda)\|_2$)  over the logarithmic penalty term
($- \sum_{(i,j)\in I_{nz}}  \ln(|x_{ij}|)$).
The convenient choice of $\tau>0$ will be illustrated for the variety of experiments considered in Section \ref{NumRes}.} \\

     Since each one of the considered variables in (\ref{optP}) is within a closed and bounded interval
    	 in $\R$, then  the feasible set $\Omega$ is a compact set.
        Although the bounds defining the set $\Omega$  ensure that the
   	 given thresholds for the entries in  $I_{nz}$ are not violated, it is worth noting that
   	 the  intervals $[l_{ij}, u_{ij}]$
   	 can be very large. In this sense, the second term of the objective function ($- \sum_{(i,j)\in I_{nz}}
   	  \ln(|x_{ij}|)$) is helpful in selecting the largest possible values within these intervals, provided
   	  that the given  eigenvalues obtained in the solution matrix are the given ones (i.e.,  provided that the first
   	   term of the objective function (\ref{ObjF}) is zero). Furthermore, compared to the computational cost
   	    required to evaluate the first term of $f(X)$ in (\ref{ObjF}), the cost of that second term is negligible.
   	In Section  \ref{NumRes}, we will demonstrate the importance of preserving both the bounds defined by $\Omega$
   	 and the second term of the function $f(X)$.   \\

 Concerning the smoothness of the objective function (\ref{ObjF}), some comments are in order.
The first term of $f(X)$  is a norm comparing the vector of obtained eigenvalues with the
 fixed given ones, and the second term is the sum of negative $\ln(|x_{ij}|)$,
 for the indices $(i,j) \in I_{nz}$. Each summand of the second term is the composition of two
 convex functions and so each one is convex, and so the second term  is a convex function.
 	Since the absolute value is Lipschitz continuous (but not differentiable)  and the logarithmic function is
   Lipschitz continuous on any closed and bounded interval contained in $(0,+\infty)$,
  then the second term of $f(X)$  is Lipschitz continuous although not differentiable.  \\

  The first term in (\ref{ObjF}) is continuous because it is the composition of a norm
  (Lipschitz continuous function) with the subtraction of the eigenvalues as a function of the entries
  of the matrix minus the vector of given eigenvalues. It is well-known that eigenvalues
  (or roots of a polynomial) are continuous as functions of the coefficients of that polynomial,
  i.e., as functions of the matrix entries; see, e.g., \cite{Harris87, Rahman02}. For a recent proof
  of this result that only requires basic algebra and the definition of continuity, see also \cite{NathRoss24}.
 It is also clear that polynomial roots as functions of their polynomial coefficients are non-convex functions.
 Concerning the differentiability of   $(\Lambda(X) - \Lambda)$, if all the given eigenvalues are simple (multiplicity 1)
 then it is indeed differentiable. The expression of the first derivative can be found although this is a
  rather costly process, as it requires the full spectral decomposition.  Both the proof of differentiability
  and the  expression of the derivatives are based upon the implicit theorem of calculus \cite{Magnus85}.
   Now,  for our special inverse problem, the multiplicity of some of the given
 eigenvalues could be strictly greater than 1. In such cases, the implicit theorem of calculus does not hold.
  Consequently,  $(\Lambda(X) - \Lambda)$ is not  differentiable, but fortunately in the symmetric case
  it is Lipschitz continuous \cite[Th. 2.4]{Lewis96}.  Additional related results in finite and infinite
  dimensional spaces can be found in \cite{Cox95, Lewis99} and references therein. In Remark \ref{Remk1}
   we  summarize the relevant characteristics of problem (\ref{optP}).

\begin{remark} \label{Remk1}
We are minimizing a non-convex continuous function on a compact set, and so the simple-bound constrained
 problem (\ref{optP}) has (local and global) solutions.  Moreover, regardless of the multiplicity of the given real
  eigenvalues, the objective function $f(X)$ is  Lipschitz continuous and bounded below on the feasible
   region $\Omega$.
\end{remark}

Finally, regarding the proposed optimization problem, we note that derivatives are absent, the objective
function is costly to evaluate, and the considered inverse matrix problem typically has few variables, so
DFO schemes are clearly a convenient choice for solving problem (\ref{optP}).

\section{A derivative-free DDS method for solving the optimization model} \label{dds}

  Our main objective is to demonstrate that  solving the optimization model (\ref{optP}) yields
  solutions to inverse structured symmetric matrix problems with fixed eigenvalues, as described in the
  previous sections. To achieve this, a variety of effective and adaptable DFO schemes have been proposed and
    analyzed in the literature; see, e.g.,
    \cite{AudetDennis06, AudetHare, ConnEtAl, Dzahini, Fasano, Huyer, Kolda, Larson, Lucidi02, Powell98}
    and references therein.
 Furthermore, variants have been  developed for several of them to solve non-differentiable Lipschitz
continuous problems. Additionally, these schemes typically have user-friendly software packages available;
 see, e.g., \cite{Audet22, Liu22, Muller, Powell09, Wild}.
Any of these DFO schemes, suitable for solving non-smooth problems, can be adapted to solve the optimization
 model (\ref{optP}) with greater or lesser effort. Nevertheless, for the specific matrix application we are
  addressing, the DDS family of DFO methods is a convenient choice that is well suited to the presence of
   simple-bound constraints. We are particularly interested in the GLODS variant of the DDS family because it
    possesses a convergence result in the non-differentiable yet Lipschitz-continuous case and tends to converge
     to global solutions \cite{Custodio15}.  \\

In general, a typical iteration of a DDS scheme involves two steps:  a search step and a poll step.
The search step is optional and not required to guarantee convergence.  It can involve any finite
	 strategy to generate a set of points  (e.g., a grid, a mesh, or even randomly chosen points),  in an
	  attempt to cover the whole feasible region.
If the search step fails to improve the current iterate, the more disciplined poll step is mandatory
to obtain convergence results.  At this step, in the presence of non-smooth functions, the selection
  of normalized directions asymptotically dense in the unit sphere is a crucial feature for the method to
  converge \cite{AudetDennis06, Vicente}.  \\

The poll step evaluates the points associated with the poll directions, scaled by a step size parameter,
 in an attempt to improve the function value at the current iterate. Let $X_k$ be the current iterate
and $\delta_k>0$ the current step size. The poll step evaluates the function $f$ at the points in the set
$$ P_k = \{X_k + \delta_k d : d \in {D}_k \}, $$
where ${D}_k$ denotes the set of poll directions considered at iteration $k$. This evaluation procedure can be
performed under opportunistic or complete strategies; see, e.g., \cite{ConnEtAl}.
 In the former, the polling procedure is stopped once
that a poll point that improves the function value at the current iterate is found. In the latter, all poll points
are evaluated and the best one is chosen. If an iteration is successful, meaning that a point was found, either
at the search or the poll steps, with
a better function value than the one of the current iterate, the new point is accepted as the current iterate
and the step size is kept constant or can be increased.
 At unsuccessful iterations, the step length must be reduced. Similarly to other derivative-free optimization algorithms, DDS makes use of the extreme barrier
function $f_{\Omega}(X)$ defined as
   \begin{equation} \label{ExtF}
f_{\Omega}(X) = \begin{cases}
	f(X) & \;\mbox{ if }\; X\in\Omega,\\
	+\infty & \;\mbox{ otherwise},
\end{cases}
\end{equation}
where   $\Omega$ is the compact feasible region given by (\ref{Omega}). \\

  In either the search or poll step, we say that  $X$ is better than or improves the function
	value at $X_k$ if a simple reduction is achieved, i.e., if $f_{\Omega}(X) < f_{\Omega}(X_k)$, or if
	  a more demanding reduction (also known as a sufficient decrease) is observed, which means defining a
	   continuous and non-decreasing forcing function $\bar{\rho}$ such that
	  $f_{\Omega}(X) < f_{\Omega}(X_k) - \bar{\rho}(\delta_k)$,  where $\delta_k>0$ is the current step length.
 For further details about the search step and the poll step, as well as the different versions
 of DDS methods, we recommend looking at \cite{AudetDennis, AudetHare, ConnEtAl, Dzahini, Kolda}.
 For the sake of completeness, a general DDS scheme for solving (\ref{optP}) is  presented in
  Algorithm  \ref{ddsA}.  \\

\begin{algorithm}
	\begin{rm}
		\begin{description}
			{\small \vspace{1ex}
				\item[Initialization] \ \\
			Choose a set  $\cal D$  (possibly infinite) of positive spanning sets,  and   the initial step size
				$\delta_0>0$.
				 Define $0<\epsilon_1< 1$, the
				coefficient for step size contraction.  Consider
				$X_0$ in the feasible set as the initial guess.
				
				\item[For $k=0,1,2,\ldots$] \
				\begin{enumerate}
				 \item[1.] {\bf Search step:} Try to compute a point $X\in \Omega$ such that
				 $X$ improves the value of $f_{\Omega}(X_k)$ by evaluating
				  $f_{\Omega}$ at a finite number of points (chosen in an organized way or randomly).
				  If such a point is found then set $X_{k+1}=X$,
				 declare the iteration and the search step successful, and skip the poll step. \ \\
					
					\item[2.] {\bf Poll step:} Choose a positive spanning set ${D}_k$
					from  $\cal D$. Evaluate the extreme barrier function $f_{\Omega}$ at the
					poll set $P_k=\{X_k + \delta_k d : d \in D_k\}$. If a poll point is found such that
					 $f_{\Omega}(X_k +\delta_k d)$ improves the value of $f_{\Omega}(X_k)$, then set
					  $X_{k+1}=X_k+\delta_k d$, and declare the poll step as successful. Otherwise,
					   the poll step is declared as unsuccessful.\ \\
					
					\item[3.] {\bf Step size parameter update:} If
					the poll step was unsuccessful, then reduce the step size parameter,
					$\delta_{k+1}=\epsilon_1\delta_k$, otherwise set $\delta_{k+1}=\delta_k$.\
								
				 \end{enumerate}
				\item[EndFor]\ }
		\end{description}
	\end{rm}
	\caption{: a general DDS  algorithm for solving (\ref{optP}).}
	\label{ddsA}
\end{algorithm}

As a particular case of Algorithm 1, in this work we will consider a special variant called GLODS
\cite{Custodio15}, which has certain features that allow
 it to frequently converge to global solutions. Among these features, a very clever multistart strategy is
  particularly noteworthy.  Moreover, GLODS has theoretical support for the non-smooth
  	characteristics of problem (\ref{optP}) as described in Remark \ref{Remk1}. \\

 As with any DDS scheme,  each iteration of GLODS is structured around a search step and/or a poll step. The objective is to
 thoroughly explore promising areas of the feasible region $\Omega$ using local search (the poll step), while
   identifying these areas using multistart strategies (the search step). For this purpose, the extreme barrier function $f_{\Omega}$
   described before is used.
 The main idea behind the GLODS strategy is to start multiple local search lines from different points distributed
  over the feasible region. These points are computed by several possible pattern sampling strategies.
  Each evaluated point has a comparison radius which is used to merge different local lines of search
  when they are sufficiently close to each other.
  Therefore, the new lines of search are not always conducted all the way until the end. More specifically,
   proximity of the points detected by the comparison radius prevents points classified as inactive from
    being further explored.
    The goal of GLODS is to end up with as many active points as the number of local minimizers,
     which would facilitate the
   straightforward identification of the potential global extreme value.  This combination of features enables GLODS to achieve
   greater efficiency than traditional multistart strategies. For a full description of the algorithm, and further information on the
   features of GLODS, please see \cite{Custodio15}. \\

  Concerning the convergence properties of GLODS,  we recall that, as noticed in Remark \ref{Remk1},
  the function $f(X)$ in problem (\ref{optP}) is Lipschitz continuous and bounded below. In that case, the
  convergence properties of GLODS have been analyzed  in \cite[Section 3]{Custodio15}.
  Indeed, it is  established that the Clarke generalized directional derivatives
  \[ f^o(X;\widehat{d}) \equiv \limsup_{X'\rightarrow X,\; X'\in\Omega,\; t\downarrow 0,\; X'+t\widehat{d} \in \Omega}
  \frac{f(X'+t\widehat{d}) - f(X')}{t}, \]
   computed at a limit point $X^*$ of a refining subsequence of the sequence of iterates generated by the
    algorithm, are nonnegative (i.e., $f^o(X^*;\widehat{d})\geq 0$) for
  all directions $\widehat{d}$ in the interior of the Clarke  tangent cone $T_{\Omega}^{Cl}(X^*)$ to the
   compact feasible region $\Omega$ given by (\ref{Omega}).
  Roughly speaking, $T_{\Omega}^{Cl}(X^*)$   is the set of directions that start at $X^*$ and point into the  set $\Omega$.
   For a formal definition of a refining subsequence and the Clarke tangent cone, refer to
  \cite{AudetDennis, AudetHare, Clarke, Custodio15}.   By assuming density of the set of refining
   directions in the unit sphere, the fact that $f^o(X^*;\widehat{d})\geq 0$ for all $\widehat{d}\in T_{\Omega}^{Cl}(X^*)$ establishes that there is a limit point that is Clarke stationary in the
   non-differentiable constrained case \cite{Clarke}.
     For more details on the mathematical formalities of these last statements, please see  \cite{AudetDennis, AudetDennis06, AudetHare, Clarke, Custodio15}.

\section{Heuristic strategies for solving the optimization model} \label{Heurist}

 Another type of DFO iterative methods that can be used to solve problem (\ref{optP}) are the so-called
  (meta) heuristic strategies. The central common feature of all heuristic optimization methods is that they
  start  with a more or less arbitrary feasible initial guess, iteratively produce new solutions by combining
   randomness and a generation rule, evaluate these new solutions using the objective function, and
   eventually report the best feasible solution found during the search process; see, e.g.,
    \cite{banga, Chopard, Talbi}  and references therein. The execution of the
   iterated search procedure is usually stopped when there has been no further improvement over a given
    number of iterations (or further improvements cannot be expected), or when the found solution is
    good enough. Specifically we will consider the following well-known strategies:
  Genetic Algorithms, Particle Swarm Optimization,   and Covariance Matrix Adaptation Evolution
   Strategy (CMA-ES). \\

   A Genetic Algorithm (GA) is a technique for addressing both constrained and unconstrained optimization
    challenges,  inspired by the natural selection process seen in biological evolution. The algorithm
    iteratively  updates a     population of potential solutions. In each iteration, individuals are
    randomly chosen from the existing population  to act as parents, generating offspring for the
      subsequent generation, and also by accepting random mutations.
     Through repeated generations, the population gradually progresses towards an optimal
     solution \cite{Chopard, Hansen}.   \\

 Particle Swarm Optimization (PSO) is a heuristic optimization method which is biologically-inspired to simulate
  the behavior of swarms, for example, a
flock of birds or a school of fish. The process begins with a randomly generated population of candidate solutions, which are
 referred to as particles. These particles are then subjected to iterative movement within the search space, guided by the
  application of straightforward mathematical formulas that consider both their current position and velocity. The movement of
   each particle is influenced by its local best known position; however, it is also guided towards the best known global positions
    in the search space. These global positions are updated as better positions are found by other particles. It is expected that
    this scheme will facilitate the convergence of the swarm towards optimal solutions. The parameters or weights associated with
     the local best known position and the global one can be updated in a dynamic and stochastic manner throughout the search process;
      see, e.g., \cite{Kennedy, Poli}.  PSO shares many similarities with evolutionary techniques: the iterative process is
       initialized with a population of random  solutions and it searches for optima by updating generations. However,
 unlike evolutionary schemes, PSO has no evolution operators such as crossover and mutation. \\

 	CMA-ES is a stochastic method that belongs to the general family of evolutionary algorithms for
 	 which the search strategy is not based on the use
 of function derivatives, but instead it is only based on the use of  objective function evaluations, that in our case is given by $f(X)$ in
 (\ref{ObjF}). In general, an initial population  is generated, and then at every cycle new solution candidates or offsprings are generated
 from the solutions already stored in the population matrix. The new candidates are obtained by using random biologically-inspired combinations of
 their parents (crossover),  and also by using random mutations; see, e.g., \cite{Hansen, Hansen03}.
 At every cycle, based on the objective function evaluation, the new candidates are either accepted to be part of the current population, or rejected.
 An attractive feature of  CMA-ES is that a covariance matrix is used to control the generation of new candidates. Moreover, the
 CMA-ES approach includes a dynamic self-updating process of the  required parameters which are used for the generation  of offsprings.

\section{Numerical Results} \label{NumRes}

To  illustrate the practical
performance of the proposed DFO algorithms, we present the results of some numerical experiments.
For each considered problem we will explore different bounds on some of the entries in the $I_{nz}$ index set.
All the experiments were run using Matlab R2024b with double precision, and they were executed in a laptop
 computer with CPU Intel core i7-1165G7, 4.7 GHz, and 16GB of RAM memory.

 \subsection{Performance of GLODS} \label{glods}

 We will now exhibit the performance  of  GLODS on several experiments.
  GLODS is implemented in MATLAB and freely available for use at
 	{\bf  https://docentes.fct.unl.pt/algb/pages/glods}, under a GNU Lesser General Public License.
  The runs  terminate using the  GLODS stopping criteria based on a fixed tolerance $tol>0$  {or
  when an appropriate precision of the eigenvalues is achieved.}
 The required input parameters  are fixed as follows:
 $tol= 10^{-5}$, $\delta_0=1$, $\epsilon_1 = 1/2$,  and  $\tau=2*\max(|\lambda_i|)$ (i.e., twice the largest
  in absolute value of the given eigenvalues).  In all our experiments, this value of the parameter
  $\tau$ is large enough to guarantee that the eigenvalues of the solution matrix are {the given
   ones with the required precision.  This is an empirical rule that was found to be effective in this study.}
    	To illustrate the convenience of this choice, in all the  experiments we will show that, for a sequence of
   values ranging from nearly zero to less than half the aforementioned value, GLODS obtains solutions for which
    the first term of $f(X)$ is not zero.
   In all cases, the evaluation procedure is performed under opportunistic strategies, that is, the polling
    procedure is stopped, once a poll point that improves the function value at the current iterate is found. 
  In our  experiments, the GLODS parameter dir$\_$dense is set to 1 allowing  a random rotation of the positive
   	  basis, $[I;-I]$,  at each iteration, where $I$ denotes the identity matrix. This choice in the polling step
   	   guarantees, asymptotically, a dense set of directions on the unit sphere; see \cite{Custodio15}.
    To address the box-constraints we used the extreme barrier
   approach described in Section \ref{dds}.
     GLODS is initialized with $n$ points, conveniently located in a
   line joining the lower and upper bound vectors ($l\equiv l(i,j)$ and $u\equiv u(i,j)$)  that define  $\Omega$.
   In the following iterations, by setting the GLODS parameter search$\_$option $=1$,
   	 we use a search step employing a Latin hypercube sampling strategy that generates asymptotically a dense
   	 set of points; for the details see \cite{McKay}.
   In all experiments, the algorithms are stopped if either the number of iterations or the number of
   function evaluations reached $3000$ or $3000\times n$, respectively,  or if $\delta_k$ was less than
   $tol=10^{-5}$, {or if $\|vec(\Lambda(X_k) - \Lambda)\|_2\leq tol$.}
   In each case,  we report the sum in absolute value of all the entries in $I_{nz}$ that appear
   	in the solution matrix (sum$(|NZ|)$),
   	 and also the smallest entry in $I_{nz}$ in absolute value ($\min(|NZ|)$). We also report
   the required number of function evaluations (NumEvalf) and iterations (Iter). \\

    For our  first experiment, let us consider  the $4\times 4$ symmetric matrix indicated in \cite[p. 265]{sutton}, in which the only
   nonzero fixed value is $\pi$ at the $(2,2)$ entry.
   \[ \left[ \begin{array}{rrrr}
   	{\rm x} & {\rm x} & {\rm nz} & 0  \\
   	{\rm x} & \pi  & {\rm nz} & {\rm x} \\
   	{\rm nz} & {\rm nz} & {\rm x} &  0  \\
   	0 & {\rm x} & 0 & {\rm nz} \\
   \end{array} \right]. \]
   We note that the associated undirected  graph  is not a tree  since it has a
   cycle (vertices 1, 2, and 3). In this case, we fix the given eigenvalues to be 1, 2, 3, and 4.   Since the four eigenvalues are different,
   the existence of a symmetric matrix with the described pattern is guaranteed \cite[Chapter 11]{JohnSaiago}.
    Notice that  in this case we have 4 free variables in the index set  $I_{free}$: $(1,1), (1,2), (2,4)$,
     and $(3,3)$, and 3 nonzero variables in the set  $I_{nz}$: $(1,3), (2,3)$, and $(4,4)$.
      We will explore several possible intervals for the 3 variables in $I_{nz}$. \\

    In Table \ref{Exp1}  we report the behavior of GLODS  and the options considered for this
    $4\times 4$  matrix. In all cases,  all the structural requirements are satisfied and the 4 obtained
    eigenvalues {achieve the target eigenvalue tolerance, i.e.,
    $\|vec(\Lambda(X) - \Lambda)\|_2\leq 10^{-5}$.}
    We now present (with up to 5 digits)  the different solutions  produced
     by GLODS when the intervals for the variables   are the following:
    (top-left) $-10\leq a_{ij} \leq 10$ for $(i,j) \in I_{free}$ and $0.4 \leq a_{ij} \leq 10$ for
     $(i,j) \in I_{nz}$;  (top-right)  $-5\leq a_{ij} \leq 5$ for $(i,j) \in I_{free}$ and
     $0.4 \leq a_{ij} \leq 5$ for  $(i,j) \in I_{nz}$;
     (bottom-left) $-3\leq a_{ij} \leq 3$ for $(i,j) \in I_{free}$ and
     $0.4 \leq a_{ij} \leq 3$ for  $(i,j) \in I_{nz}$;
     (bottom-right) $-5\leq a_{ij} \leq 5$ for $(i,j) \in I_{free}$, $-5 \leq a_{13},a_{23} \leq -0.4$,  and
     $0.4 \leq a_{44} \leq 5$:
       {\small \[ \left[ \begin{array}{rrrr}
     		2.9208 &   -0.905 &   0.6595 &  0 \\
     		-0.905   & 3.1416 &   0.7034 &  -0.410 \\
     		0.6595  &  0.7034 &   1.8878 &  0 \\
     		0  &  -0.410 &  0 &   2.0497 \\
     	\end{array} \right], \;\;\;\;
     	\left[ \begin{array}{rrrr}
     	   1.8707 &  -0.1959  &  0.9489   &    0 \\
        	-0.1959  &  3.1416  &  0.8012  &  0.7052 \\
     	0.9489  &  0.8012 &   2.6286 &    0 \\
     		   0  &  0.7052 &        0 &   2.3591  \\
     	\end{array} \right]   \] }
     {\small \[ \left[ \begin{array}{rrrr}
     	2.1009 &  0.625 & 0.8818 &  0 \\
     		0.625  &  3.1416  & 0.8024 &  0.0384 \\
     	0.8818 &  0.8024 &   1.7573 &  0 \\
     		0  & 0.0384 &   0  &  3.00 \\
     	\end{array} \right], \;\;\;\;
     	\left[ \begin{array}{rrrr}
     	 1.7387   &  0.2544  &   -0.7796  &  0 \\
     	0.2544 &   3.1416  & -0.9488  &   0.4196 \\
     	 -0.7796  &  -0.9488 & 2.0579  &  0 \\
     		0  & 0.4196  &  0 &   3.0617 \\
     	\end{array} \right].   \] }

  \begin{table}[H]
  	\caption{Performance of GLODS for the $4\times 4$  symmetric matrix and all the considered bounds.}  \label{Exp1}
 	\begin{center}
  		\begin{tabular}{ccccc} \hline
  			Bounds  &  Iter &  sum$(|NZ|)$  &  NumEvalf & $\min(|NZ|)$  \\ \hline
  		$-10\leq I_{free} \leq 10$, $\:0.4 \leq I_{nz} \leq 10$ & 105 &  $4.775$ & 1765 &  0.6595 \\
  		$-5\leq I_{free} \leq 5$, $\:0.4 \leq I_{nz} \leq 5$ & 108  &  $5.859$ & 1531  & 0.8012 \\
  		$-3\leq I_{free} \leq 3$, $\:0.4 \leq I_{nz} \leq 3$ &  112 &  $6.3684$ & 2318  & 0.8024  \\
  		$-5\leq I_{free} \leq 5$, $\:-5 \leq a_{13},a_{23} \leq -0.4$ &  89 &  $6.5185$ & 1625  & 0.7796  \\
  		and $\:0.4 \leq a_{44} \leq 5$ & & & & \\
  		\hline
  		\end{tabular}
  	\end{center}
  \end{table}

 Now, fixing the   bounds as follows:  $-5\leq a_{ij} \leq 5$ for $(i,j) \in I_{free}$ and
 $0.4 \leq a_{ij} \leq 5$ for  $(i,j) \in I_{nz}$, in  Table \ref{Exp1tau} we report the norm of the
  difference between the obtained eigenvalues and the given matrix $\Lambda$
  (i.e., we report $\|vec(\Lambda(X) - \Lambda)\|_2$), {and the absolute value of the
  	 smallest entry in $I_{nz}$ ($\min(|NZ|)$)},  for an increasing set of values of the parameter
   $\tau>0$ that appears in (\ref{optP}). For small values of $\tau$, we observe that the obtained matrix
   fails to reproduce the given eigenvalues {with the required precision}, which is unacceptable
   	for the inverse problem under consideration.    {Moreover, when $\tau$ is very close to zero, the first term of the objective function (\ref{ObjF}) becomes
  practically irrelevant. In this case, the minimization mainly considers  the second term in (\ref{ObjF}), which
   is associated with the set $I_{nz}$. Hence, the obtained value of $\min(|NZ|)$ is  at the boundary
    of the feasible box $\Omega$. As we increase the value of $\tau$,   the first term in 
     (\ref{ObjF}) becomes more relevant.  Consequently, the obtained eigenvalues get closer to the given ones, 
    and the value of $\min(|NZ|)$ decreases.} 
    Finally, as reported in Table \ref{Exp1}, for $\tau = 2(\max(|\lambda_i|))$, the obtained matrix reproduces
     the given eigenvalues {up to the required tolerance}.

   \begin{table}[H]
  	\caption{Performance of GLODS for the $4\times 4$  symmetric matrix, using the bounds $-5\leq I_{free} \leq 5$,
  		 $\:0.4 \leq I_{nz} \leq 5$, and  several values of the parameter $\tau$.}  \label{Exp1tau}
   	\begin{center}
  		 {\begin{tabular}{ccc} \hline
  				$\tau$  &  $\|vec(\Lambda(X) - \Lambda)\|_2$ & $\min(|NZ|)$  \\ \hline
                $\max(|\lambda_i|)/20$ &  $8.170$ & $4.999$ \\
  				$\max(|\lambda_i|)/10$ &  $3.637$ & $2.856$ \\
  				$\max(|\lambda_i|)/5$ &  $0.974$ & $1.385$ \\
  				$\max(|\lambda_i|)/2$ & $4.0\times 10^{-5}$ & $0.9879$ \\
  				$\max(|\lambda_i|)$ & $8.2\times 10^{-6}$ & $0.8109$ \\
  		    	$2(\max(|\lambda_i|))$ & $4.1\times 10^{-6}$ & $0.8012$  \\
  				\hline
  			\end{tabular}}
  	\end{center}
  \end{table}

 As we discussed in Section 2, preserving both the simple bounds that define $\Omega$ and the second
 term of the function $f(X)$, i.e., $- \sum_{(i,j)\in I_{nz}}  \ln(|x_{ij}|)$, is valuable. To illustrate
 this, we applied GLODS to the four cases reported in Table \ref{Exp1}, considering that the objective
 function is given by $f(X) = \|vec(\Lambda(X) - \Lambda)\|_2$. The results show that, by maintaining
 essentially the same number of iterations and  function evaluations, the average of sum$(|NZ|)$
 is reduced from 6.08 to 4.53 and $\min(|NZ|)$  is reduced from 0.74 to 0.62.  At this point, we would like
  to emphasize  that the cost of the second term of $f(X)$ in (\ref{ObjF}) is negligible compared to the cost of
   evaluating $\|vec(\Lambda(X) - \Lambda)\|_2$. \\

   Additionally, it is worth mentioning that, for this $4\times 4$
    experiment, the optimization approach proposed in \cite{chehab} reports matrices that solve the same problem
     (Experiment 3 in \cite[pp. 20-21]{chehab}). However,
    since they cannot impose bounds on their manifold optimization schemes, they obtain an average sum$(|NZ|)$
     of approximately 4.8 and a $\min(|NZ|)$ of approximately 0.2. In some cases, the obtained $\min(|NZ|)$ is as
     low as 0.1. These results are lower than those reported by our DFO approach, particularly the $\min(|NZ|)$
     value, which is the most significant. The ability to obtain solution graphs (or matrices) with higher
      edge weights, in the $I_{nz}$ set, that  always exceed the minimum required bounds is clearly an advantage
       of the proposed DFO-based approach. \\

For the second experiment, we consider the $7\times 7$ symmetric Laplacian matrix described in
 \cite[Table II, p. 369]{Trinaj}. These matrices have a potential for use in chemical graph theory.
  We fix the diagonal entries of $\Lambda$ to be the eigenvalues of the  matrix displayed in \cite{Trinaj}, i.e.,
   the given eigenvalues are: 0, 0.3672, 1.0571, 2.3745, 4.4681, 7.0357, and 8.6973.
    The required structural pattern is the following:
\[ \left[ \begin{array}{rrrrrrr}
	{\rm x} & {\rm nz} & 0 & 0 & 0 & 0 & 0 \\
	{\rm nz} & {\rm x} & {\rm nz} & 0 & 0 & 0 & -3 \\
	0 & {\rm nz} & {\rm x} & {\rm nz} & 0 & -3 & 0  \\
	0 & 0 & {\rm nz} & {\rm x} & {\rm nz} & 0  & 0 \\
	0 & 0 & 0 & {\rm nz} & {\rm x} & 0 & 0  \\
	0 & 0 & -3 & 0 &   0 &  {\rm x} & 0 \\
	0 & -3 & 0 & 0 &   0 &  0 & {\rm x} \\
\end{array} \right], \]
which is identified with the  undirected  weighted graph (a linear tree in this case) of the 2,3-Dimethyl pentane, shown in
Figure \ref{fig:2,3-Dimethyl pentane}.
In Table \ref{Exp2}  we report the behavior  of GLODS and the considered options for this
 $7\times 7$  Laplacian matrix. In all cases,
the structural requirements are satisfied and the 7 eigenvalues obtained are the given ones
 {with the required precision}.
  We now present (with up to 5 digits)  four different solutions  produced by GLODS, which are reported in
  Table \ref{Exp2} in the order of the rows: top left, top right, bottom left, and bottom right.

 {\footnotesize \[ \left[ \begin{array}{rrrrrrr}
 		1.984  &   1.384 &   0 &  0&0&0&0 \\
 		1.384   & 3.996 &   0.5411 &  0 &0&0&-3\\
 		0  &  0.5411 &   5.252 &  2.506 &0&-3&0\\
 		0  &  0 &  2.506 &   4.827 & 1.280 &0&0\\
 		0&0&0& 1.280 &1.066 &0&0\\
 		0&0&-3&0&0&3.1267 & 0\\
 		0&-3&0&0&0&0&3.748\\
 	\end{array} \right], \;\;\;\;
 	\left[ \begin{array}{rrrrrrr}
    	2.259 &  1.758 &  0 &  0 &0&0&0\\
 		1.758  &  5.786 &   1.158 &  0&0&0&-3 \\
 		0  & 1.158 &  5.578 &   0.9515 &0&-3&0\\
 		0  &   0 &   0.9515 &   2.893 &2.038&0&0 \\
 		0&0&0&2.038&2.0583&0&0\\
 		0&0&-3&0&0&2.7534&0\\
 		0&-3&0&0&0&0&2.6725
 	\end{array} \right]   \] }
 {\footnotesize \[ \left[ \begin{array}{rrrrrrr}
 		0.573 & -0.763 & 0 &  0 &0&0&0\\
    	-0.763  &  5.157  & -0.823 &  0 &0&0&-3\\
 		0 &  -0.823 &   5.9107 &  -2.225 &0&-3&0\\
 		0  & 0 &   -2.225  &  2.743&-1.120&0&0 \\
 		0&0&0&-1.120&3.945&0&0\\
 		0&0&-3&0&0&2.97&0\\
 		0&-3&0&0&0&0&2.7
 	\end{array} \right], \;\;
 	\left[ \begin{array}{rrrrrrr}
 		2.397   &  -1.286  &   0  &  0 &0&0&0\\
 		-1.286 &   5.655  & -1.112  &   0&0&0&-3 \\
 		0  &  -1.112 &  6.1 &  -1.116 &0&-3&0\\
 		0  & 0  &  -1.116 &  2.85&2.033&0&0 \\
 		0&0&0&2.033&2.198&0&0\\
 		0&0&-3&0&0&2.55&0\\
 		0&-3&0&0&0&0&2.248
 	\end{array} \right].   \] }

 \begin{table}[H]
 	\caption{Performance of GLODS for the $7\times 7$  Laplacian matrix and all the considered bounds.}  \label{Exp2}
 	\begin{center}
 		\begin{tabular}{ccccc} \hline
 			Bounds  &  Iter &  sum$(|NZ|)$  &  NumEvalf & $\min(|NZ|)$  \\ \hline
 			$-10\leq I_{free} \leq 10$, $\:0.4 \leq I_{nz} \leq 10$  & 277 &  $11.42$ & 3362 &  0.5411 \\
 			$0.2\leq I_{free} \leq 7$, $\:0.7 \leq I_{nz} \leq 7$ &  584 &  $11.81$ & 7989  & 0.9515  \\
 			$0.2\leq I_{free} \leq 7$, $\:-7 \leq I_{nz} \leq -0.5$ &  650 &  $9.862$ & 9292  & 0.763  \\
 			$1\leq I_{free} \leq 7$, $\:-7 \leq I_{nz} \leq -1$ &  385 &  $11.094$ & 5222  & 1.112  \\
 			\hline
 		\end{tabular}
 	\end{center}
 \end{table}

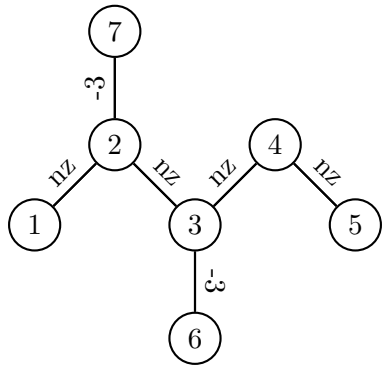
\begin{figure}[H]
    \centering
\begin{tikzpicture}[node distance={15mm}, thick, main/.style = {draw, circle}]
\node[main] (1) {$1$};
\node[main] (2) [above right of=1] {$2$};
\node[main] (3) [below right of=2] {$3$};
\node[main] (4) [above right of=3] {$4$};
\node[main] (5) [below right of=4] {$5$};
\node[main] (6) [below of=3] {$6$};
\node[main] (7) [above of=2] {$7$};
\draw (1) -- (2);
\draw (2) -- (3);
\draw (3) -- (4);
\draw (4) -- (5);
\draw (3) -- (6);
\draw (2) -- (7);
\draw (1) -- node[midway, above, sloped] {nz} (2);
\draw (2) -- node[midway, above, sloped] {nz} (3);
\draw (3) -- node[midway, above, sloped] {nz} (4);
\draw (4) -- node[midway, above, sloped] {nz} (5);
\draw (2) -- node[midway, above, sloped] {-3} (7);
\draw (3) -- node[midway, above, sloped] {-3} (6);
\end{tikzpicture}
\caption{Structural pattern graph of the symmetric Laplacian matrix (2,3-Dimethyl pentane).}
    \label{fig:2,3-Dimethyl pentane}
\end{figure}

	Now, fixing the   bounds as follows:  $-10\leq a_{ij} \leq 10$ for $(i,j) \in I_{free}$ and
	$0.4 \leq a_{ij} \leq 10$ for  $(i,j) \in I_{nz}$, in  Table \ref{Exp2tau} we report the value of
	 $\|vec(\Lambda(X) - \Lambda)\|_2$  for an increasing set of values of the parameter
	$\tau>0$. As before, for small values of $\tau$ we observe that the obtained matrix
	fails to reproduce the given eigenvalues, and also  that as  $\tau$ increases the obtained
	eigenvalues get closer to the given ones. As reported in Table  \ref{Exp2}, for $\tau = 2(\max(|\lambda_i|))$,
	the obtained matrices reproduce the given eigenvalues {with the required precision}.

\begin{table}[H]
	\caption{Performance of GLODS for the $7\times 7$  Laplacian matrix, using the bounds $-10\leq I_{free} \leq 10$,
		$\:0.4 \leq I_{nz} \leq 10$, and  several values of the parameter $\tau$.}  \label{Exp2tau}
 	\begin{center}
			\begin{tabular}{cc} \hline
				$\tau$  &  $\|vec(\Lambda(X) - \Lambda)\|_2$  \\ \hline
				$\max(|\lambda_i|)/10$ &  $1.582$   \\
				$\max(|\lambda_i|)/5$ &  $0.24$ \\
				$\max(|\lambda_i|)/2$ & $0.0025$  \\
				$\max(|\lambda_i|)$ & $3.2\times 10^{-4}$  \\
				\hline
			\end{tabular}
	\end{center}
\end{table}

	For this $7\times 7$ experiment, the optimization approach proposed in \cite{chehab} reports
	matrices that solve the same problem (Experiment 2 in \cite[pp. 19-20]{chehab}). However, they obtain an
	average sum$(|NZ|)$ of approximately 8.8 and a $\min(|NZ|)$ of approximately 0.21. In some cases,
	the obtained $\min(|NZ|)$ is as low as 0.12. These results are lower than those reported by our
	DFO approach, especially  the $\min(|NZ|)$ value, which is the most important. Furthermore,
   to illustrate the value of preserving the bounds that define $\Omega$ and the second
	term of the function $f(X)$, i.e., $- \sum_{(i,j)\in I_{nz}}  \ln(|x_{ij}|)$, we applied GLODS to the four cases reported in Table \ref{Exp2}, considering that the objective
	function is given by $f(X) = \|vec(\Lambda(X) - \Lambda)\|_2$. The results show that, by maintaining
	essentially the same number of iterations and  function evaluations, the average of sum$(|NZ|)$
	is reduced from 11.05 to 8.92, and $\min(|NZ|)$  is reduced from 0.84 to 0.65. We note that
	 the tendency to produce, on average, a reduction in both sum$(|NZ|)$ and $\min(|NZ|)$
	  when considering $f(X) = \|vec(\Lambda(X) - \Lambda)\|_2$ is maintained throughout the remaining
	   experiments, as was observed in the first two. Therefore, for reasons of space, we will not report these
	    results for future experiments. \\

For the third experiment, we consider the $10\times 10$ symmetric matrix described in \cite[p. 100]{JohnSaiago}, which is identified
with the undirected graph shown in Figure \ref{fig:p150JohnSaiago}. This graph is the smallest possible nonlinear tree;
see \cite[Section 5]{JohnSaiago}. Following the description in \cite[p. 127]{JohnSaiago}, we fix the given eigenvalues to
 be -3, -2, -2, 0, 0, 0, 0, 2, 2, and 3. These 5 different eigenvalues together with their multiplicity guarantee that there
  are possible solutions to the inverse structural matrix problem; see \cite[Section 6]{JohnSaiago}.
   The required structural pattern is the following:
\[ \left[ \begin{array}{rrrrrrrrrr}
	0&0 & a & 0 & 0 & 0 & 0 & 0& 0 & 0 \\
	0 & 0 & \sqrt{2} & 0 & 0 & 0 & 0& 0 & 0& 0 \\
	a & \sqrt{2} &0&1& 0 & 0 & 0 & 0& 0 & 0  \\
	0 & 0 &1 & 0 & \sqrt{2} & 0  & 0&\sqrt{2}&0&0 \\
	0 & 0 & 0 & \sqrt{2} &0& b & c & 0 & 0 &0 \\
 	0 & 0 & 0&0& b &0&0&0&0&0\\
 	0 & 0 & 0&0& c &0&0&0&0&0\\
 	0 & 0 &0&\sqrt{2}&0&0&0&0& d & e \\
 	0 & 0 &0&0&0&0&0& d &0&0\\
 	0 & 0 &0&0&0&0&0& e &0&0\\
\end{array} \right], \]
where $a, b, c, d$, and $e$ are the only considered variables, and they must belong to the nonzero  set of indices
$I_{nz}$.

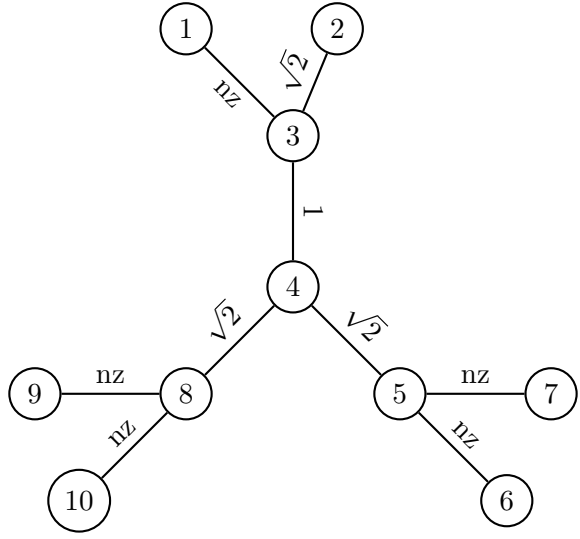
\begin{figure}[H]
    \centering
\begin{tikzpicture}[node distance={20mm}, thick, main/.style = {draw, circle}]
\node[main] (1) {$1$};
\node[main] (2) [right of=1] {$2$};
\node[main] (3) [below right of=1] {$3$};
\node[main] (4) [below  of=3] {$4$};
\node[main] (5) [below right of=4] {$5$};
\node[main] (6) [below right of=5] {$6$};
\node[main] (7) [ right of=5] {$7$};
\node[main] (8) [below left of=4] {$8$};
\node[main] (9) [left of=8] {$9$};
\node[main] (10) [below left of=8] {$10$};
\draw (2) -- node[midway, above, sloped] {$\sqrt{2}$} (3);
\draw (1) -- node[midway, below, sloped] {{\rm nz}} (3);
\draw (3) -- node[midway, above, sloped] {1} (4);
\draw (4) -- node[midway, above, sloped] {$\sqrt{2}$} (5);
\draw (5) -- node[midway, above, sloped] {{\rm nz}} (7);
\draw (5) -- node[midway, above, sloped] {{\rm nz}} (6);
\draw (8) -- node[midway, above, sloped] {{\rm nz}} (9);
\draw (8) -- node[midway, above, sloped] {{\rm nz}} (10);
\draw (4) -- node[midway, above, sloped] {$\sqrt{2}$} (8);

\end{tikzpicture}
\caption{Structural pattern graph of the $10\times 10$ real symmetric matrix described in \cite[p. 127]{JohnSaiago}.}
    \label{fig:p150JohnSaiago}
\end{figure}

 We  now present (with up to 5 digits) the four different matrices obtained by GLODS for different bounds on the nonzero
 variables $a,b,c,d$, and  $e$: \\ [2mm]
Exp3-Matrix-1 is completed with $(a,b,c,d,e)=(1.4142, 1.2509, 1.5606, 1.6640, 1.1096)$, \\  [2mm]
Exp3-Matrix-2  with $(a,b,c,d,e)=(1.4142,  1.502,  1.3206, 1.6587, 1.1174)$,  \\ [2mm]
Exp3-Matrix-3  with $(a,b,c,d,e)=(1.4142, 1.4152,  1.4133, 1.5082,  1.3135)$,  \\ [2mm]
 Exp3-Matrix-4  with $(a,b,c,d,e)=(1.4142, 1.4172, 1.4112, 1.2526,  1.5592)$. \\ [2mm]

 \begin{table}[H]
 	\caption{Performance of GLODS for the  symmetric matrix in \cite[p. 127]{JohnSaiago} and all the considered bounds.}  \label{Exp3}
 	\begin{center}
 		\begin{tabular}{cccccc} \hline
 		Solution &	Bounds  &  Iter &  sum$(|NZ|)$   &  NumEvalf & $\min(|NZ|)$  \\ \hline
 	Exp3-Matrix-1 &	$0.5 \leq I_{nz} \leq 10$   & 124 &  $13.998$ & 776 & 1.1096 \\
    Exp3-Matrix-2 & $0.5 \leq I_{nz} \leq 5$ &  100 &  $14.026$ & 662  & 1.1174  \\
 	Exp3-Matrix-3 & $1 \leq I_{nz} \leq 5$ & 106 &  $14.13$ & 676  & 1.3135  \\
 	Exp3-Matrix-4 &	$1 \leq I_{nz} \leq 3$  & 107 &  $14.08$ & 733  & 1.2526  \\
 				
    \hline
 		\end{tabular}
 	\end{center}
 \end{table}

	Now, fixing the   bounds as follows:  $1 \leq a_{ij} \leq 5$ for  $(i,j) \in I_{nz}$, in
	Table \ref{Exp3tau} we report the value of
	$\|vec(\Lambda(X) - \Lambda)\|_2$  for an increasing set of values of the parameter
	$\tau>0$. It can be observed that for small values of $\tau$ the obtained matrix
	fails to reproduce the given eigenvalues {with the required precision}, and also  that
	as  $\tau$ increases the obtained
	eigenvalues get closer to the given ones. As reported in Table  \ref{Exp3}, for $\tau = 2(\max(|\lambda_i|))$,
	the obtained matrices reproduce the given eigenvalues {with the required tolerance}.

\begin{table}[H]
	\caption{Performance of GLODS for the symmetric matrix in \cite[p. 127]{JohnSaiago},
		using the bounds $1 \leq I_{nz} \leq 5$, and  several values of the parameter $\tau$.}  \label{Exp3tau}
	\begin{center}
			\begin{tabular}{cc} \hline
				$\tau$  &  $\|vec(\Lambda(X) - \Lambda)\|_2$  \\ \hline
				$\max(|\lambda_i|)/10$ &  $10.53$   \\
				$\max(|\lambda_i|)/5$ &  $4.006$ \\
				$\max(|\lambda_i|)/2$ & $1.2\times 10^{-5}$  \\
				$\max(|\lambda_i|)$ & $2.6\times 10^{-6}$ \\
				\hline
			\end{tabular}
	\end{center}
\end{table}

	For this $10\times 10$ real symmetric matrix, the optimization approach proposed in \cite{chehab} also
	reports matrices that solve the same problem (Experiment 4 in \cite[pp. 21-22]{chehab}). However, they obtain
	 an average sum$(|NZ|)$ of approximately 12.8 and a $\min(|NZ|)$ of approximately 0.2. These results are
  lower than those reported by our DFO approach as observed in Table \ref{Exp3}, especially  the $\min(|NZ|)$
   value, which is the most significant. \\

For the fourth experiment, we  fix the diagonal entries of $\Lambda$ to be the eigenvalues of the tridiagonal symmetric positive definite $n\times n$
matrix obtained by discretizing the second-order derivative, over a closed interval in one variable, with zero values at the extreme points.
In other words, the given $\lambda_i$s are the eigenvalues of the matrix
$tridiag(-1,2,-1)$,
which can be seen as a special case of the so-called Laplacian
matrices that have an impact in chemical graph theory \cite{Trinaj}.
It is well-known that, for any $n \geq 2$, the eigenvalues of $tridiag(-1,2,-1)$
are given by $\lambda_i =2-2 \cos(i \pi/(n + 1))$
for $1 \leq i \leq n$. The considered
structure is the following: outside of the three main diagonals
all entries are zero, the entries in the
main diagonal are all totally free variables (indicated with the symbol $x$), there are no nonzero preset
values, and all the entries in the main sub-diagonals are free
variables that must be nonzero (indicated with the symbol nz). For example, the matrix pattern for $n = 6$ is given by:

\[ \left[ \begin{array}{rrrrrr}
	x& {\rm nz} & 0 & 0 & 0 & 0 \\
	{\rm nz} & x &{\rm nz} & 0 & 0 & 0 \\
	0 &{\rm nz} &x&{\rm nz}& 0 & 0  \\
	0 & 0 &{\rm nz} & x &{\rm nz} & 0\\
	0 & 0 & 0 & {\rm nz} &x& {\rm nz}\\
 	0 & 0 & 0&0& {\rm nz}&x\\
\end{array} \right], \]
which is identified with a simple undirected graph of 6 vertices \cite[pp. 5-8]{JohnSaiago}. For this tridiagonal matrix
and for several values of $n$, in Table \ref{Exp4} we report the performance of GLODS.
We note that the eigenvalues of the obtained matrices are the same as the eigenvalues of
 $tridiag(-1, 2, -1)$ {up to the required precision}.

 \begin{table}[H]
 	\caption{Performance of GLODS for  the inverse structured tridiagonal matrix case and all the considered bounds.}  \label{Exp4}
  	\begin{center}
 		\begin{tabular}{c|ccccc} \hline
 			$n$ & Bounds  &  Iter &  sum$(|NZ|)$   &  NumEvalf & $\min(|NZ|)$  \\ \hline
 		&	$0.2\leq I_{free} \leq 10$, $\:0.2 \leq I_{nz} \leq 10$  & 169 &  $5.85$ & 1407 & 0.8008 \\
  4	&		$0.5\leq I_{free} \leq 5$, $\:0.5 \leq I_{nz} \leq 5$ &  227 &  $5.97$ & 1902  & 0.9773  \\
 	&	$1\leq I_{free} \leq 5$, $0.5\leq a_{34} \leq 5$, $-5 \leq a_{12}, a_{23} \leq -0.5$  &  188 &  $5.95$ & 1481  & 0.9099  \\  \hdashline
     		&	$0.2\leq I_{free} \leq 10$, $\:0.2 \leq I_{nz} \leq 10$  & 328 &  $11.79$ & 5185 & 0.832 \\
  7	&	$0.5\leq I_{free} \leq 5$, $\:0.5 \leq I_{nz} \leq 5$ &  273 &  $11.88$ & 4382  & 0.8927  \\
 	&	$0.8\leq I_{free} \leq 4$, $\:0.8 \leq I_{nz} \leq 4$ &  331 &  $11.96$ & 4353  & 0.9747  \\  \hdashline
     		&$0.2\leq I_{free} \leq 10$, $\:0.2 \leq I_{nz} \leq 10$   & 532 &  $17.97$ & 9804 & 0.9525 \\
  10	&	$0.5\leq I_{free} \leq 5$, $\:0.5 \leq I_{nz} \leq 5$ &  455 &  $17.93$ & 9108 & 0.9268  \\  %
 	&		$0.5\leq I_{free} \leq 5$,  $\:-5 \leq a_{12}, a_{56}, a_{9,10} \leq -0.5$ &  516 &  $17.94$ &  9780  & 0.877 \\
 	&	and $\:0.5\leq a_{23}, a_{34}, a_{45}, a_{67}, a_{78}, a_{89} \leq 5$,	 &   &  &   & \\
  \hline
 		\end{tabular}
 	\end{center}
 \end{table}

  For the sake of completeness, we now present (with up to 5 digits) the obtained matrices produced by GLODS,
 reported in Table~\ref{Exp4},  for the case in which the bounds are given by
 $0.5\leq I_{free} \leq 5$, $\:0.5 \leq I_{nz} \leq 5$.

For $n=4$ we obtain:
\[ \left[ \begin{array}{rrrr}
	2.17 &   1.0032 &   0&  0 \\
	1.0032  & 1.9642 &   1.0047 &  0 \\
	0  &  1.0047 & 2.0287 & 0.9773 \\
	0  &  0 &  0.9773 &  1.837 \\
\end{array} \right], \]
\vspace{2mm}
for $n=7$ we obtain:
{\small \[ \left[ \begin{array}{rrrrrrr}
		1.9179 &   1.0524 &   0&  0 &0&0&0\\
		1.0524  & 2.2312 &  0.9050 &  0&0&0 &0\\
		0  &  0.9050 & 1.8332 & 1.0877&0&0 &0\\
		0  &  0 &  1.0877 &  2.0652 &0.9094&0&0\\
		0  &  0 &  0 &  0.9094 &2.081&1.0954&0\\
		0  &  0 &  0 & 0 &1.0954&1.8206&0.8927\\
		0  &  0 &  0 & 0 &0&0.8927&2.05\\
	\end{array} \right],  \] }
\noindent
            and for $n=10$ we obtain:
                 {\footnotesize \[ \left[ \begin{array}{rrrrrrrrrr}
     		   1.9557 &  1.001  &      0    &    0&        0 &       0  &      0    &    0  &      0&        0\\
   1.001  &2.0967  & 0.9408 &       0    &    0 &       0      &  0 &       0       & 0 &       0\\
        0  &0.9408 &  1.9721  & 0.9341   &     0  &      0     &   0  &      0       & 0  &      0\\
        0   &     0  &0.9341  & 2.1318 & 1.0683   &     0    &    0   &     0      &  0   &     0\\
        0    &    0     &   0  & 1.0683  & 1.8448  &1.004   &     0    &    0     &   0    &    0\\
        0     &   0    &    0   &     0   &1.004   &2.1245  & 0.9268    &    0    &    0     &   0\\
        0      &  0   &     0    &    0    &    0  & 0.9268  & 1.9344  & 1.0281   &     0      &  0\\
        0      &  0  &      0     &   0   &     0   &     0  & 1.0281  & 2.0547  & 1.0571       & 0\\
        0      &  0 &       0      &  0  &      0    &    0  &      0  & 1.0571  & 1.8833  & 1.0063\\
        0      &  0&        0       & 0 &       0     &   0 &       0   &     0 &  1.0063  & 2.0019\\
            	\end{array} \right].   \] }

	Fixing the   bounds as follows:  $0.5\leq a_{ij} \leq 5$ for $(i,j) \in I_{free}$ and
	$0.5 \leq a_{ij} \leq 5$ for  $(i,j) \in I_{nz}$, in  Table \ref{Exp4tau} we report the value of
	$\|vec(\Lambda(X) - \Lambda)\|_2$  for an increasing set of values of the parameter
	$\tau>0$. We observe that for small values of $\tau$ the obtained matrix
	fails to reproduce the given eigenvalues {with the required precision}, and also  that
	 as  $\tau$ increases the obtained
	eigenvalues get closer to the given ones. As reported in Table  \ref{Exp4}, for $\tau = 2(\max(|\lambda_i|))$,
	the obtained matrices reproduce the given eigenvalues {with the required tolerance}.

\begin{table}[H]
	\caption{Performance of GLODS for the structured tridiagonal matrix for $n=4,7,10$,
		using the bounds $0.5\leq I_{free} \leq 5$, $\:0.5 \leq I_{nz} \leq 5$, and  several values of the parameter $\tau$.}  \label{Exp4tau}
	\begin{center}
			\begin{tabular}{cccc}
				 & $n=4$ & $n=7$ & $n=10$ \\ \hline
				$\tau$  &  $\|vec(\Lambda(X) - \Lambda)\|_2$ &  $\|vec(\Lambda(X) - \Lambda)\|_2$  &  $\|vec(\Lambda(X) - \Lambda)\|_2$  \\ \hline
				$\max(|\lambda_i|)/10$ &  $5.842$ &  $12.13$ &  $16.97$  \\
				$\max(|\lambda_i|)/5$ &  $1.696$ &  $ 4.332$  &  $7.24$  \\
				$\max(|\lambda_i|)/2$ & $3.7\times 10^{-5}$  & $7.6\times 10^{-5}$  & $3.1\times 10^{-1}$ \\
				$\max(|\lambda_i|)$ & $2.9\times 10^{-6}$ & $3.5\times 10^{-6}$ & $7.1\times 10^{-5}$ \\
				\hline
			\end{tabular}
	\end{center}
\end{table}

For our fifth experiment, we consider the problem of finding a real symmetric matrix with eigenvalues -2, -1, -1, 0, 1, 1, 2,
 and the following zero-nonzero pattern:

\[ \left[ \begin{array}{ccccccc}
	1 &  0 & a & 0 & 0 & 0 & 0  \\
	0 &  1 & a & 0 & 0 & 0 & 0  \\
	a &  a & b & c & 0 & 0 & 0  \\
	0 &  0 & c & 0 & c & 0 & 0  \\
	0 &  0 & 0 & c & -b & a & a  \\
	0 &  0 & 0 & 0 & a & -1 & 0  \\
	0 &  0 & 0 & 0 & a & 0 & -1  \\
\end{array} \right], \]
where  $a, b$ and $c$ are the only considered variables, and they must belong to the nonzero  set of
 indices $I_{nz}$.	
 We now present (with up to 5 digits) four of the obtained matrices produced by GLODS for different bounds:   \\ [2mm]
 Exp5-Matrix-1 is completed with $(a,b,c)=(0.5195, 0.5686, 1.1396)$, \\  [2mm]
 Exp5-Matrix-2  with $(a,b,c)=(0.7131, 0.408, 0.9482)$,  \\ [2mm]
  Exp5-Matrix-3  with $(a,b,c)=(0.5314, 0.5606, 1.1306)$,  \\ [2mm]
Exp5-Matrix-4  with $(a,b,c)=(-0.5191, 0.5688, -1.1399)$. \\ [2mm]
  We note that the matrices are different,
but they all have the  given  eigenvalues {with the required precision} and satisfy the
structural requirements.
 In  Table \ref{Exp5}  we report the behavior of GLODS to produce the 4 different solutions described above.

\begin{table}[H]
	\caption{Performance of GLODS for the $7\times 7$  symmetric matrix of the fifth experiment,
		 and  the considered bounds. }  \label{Exp5}
	\begin{center}
		\begin{tabular}{cccccc} \hline
	Solution &		Bounds  &  Iter &  sum$(|NZ|)$  &  NumEvalf & $\min(|NZ|)$  \\ \hline
	Exp5-Matrix-1 &		$0.2 \leq a,b,c\leq 4$  & 117 &  $9.851$ & 527 &  0.5195 \\
	Exp5-Matrix-2 &		$0.2 \leq a,b,c\leq 1$  & 119 &  $10.313$ & 571 &  0.408 \\
	Exp5-Matrix-3 &		$0.4\leq a,b,c\leq 2.5$ & 133 &  $9.895$ & 575 &  0.5314 \\
	Exp5-Matrix-4 &		$0.2 \leq b\leq 4$, $\:-4 \leq a,c\leq -0.2$ &  202 &  $9.85$ & 863  & 0.5191  \\
	\hline
		\end{tabular}
	\end{center}
\end{table}

	Fixing the   bounds as follows:  $0.2 \leq a_{ij} \leq 4$ for  $(i,j) \in I_{nz}$, in
	Table \ref{Exp5tau} we report the value of
	$\|vec(\Lambda(X) - \Lambda)\|_2$  for an increasing set of values of the parameter
	$\tau>0$. It can be observed that for small values of $\tau$ the obtained matrix
	fails to reproduce the given eigenvalues {with the required precision}, and also
	 that as  $\tau$ increases the obtained
	eigenvalues get closer to the given ones. As reported in Table  \ref{Exp5}, for $\tau = 2(\max(|\lambda_i|))$,
	the obtained matrices reproduce the given eigenvalues {with the required tolerance}.

\begin{table}[H]
	\caption{Performance of GLODS for $7\times 7$  symmetric matrix of the fifth experiment,
		using the bounds $0.2 \leq I_{nz} \leq 4$, and  several values of the parameter $\tau$.}  \label{Exp5tau}
 	\begin{center}
			\begin{tabular}{cc} \hline
				$\tau$  &  $\|vec(\Lambda(X) - \Lambda)\|_2$  \\ \hline
				$\max(|\lambda_i|)/10$ &  $9.46$   \\
				$\max(|\lambda_i|)/5$ &  $4.736$ \\
				$\max(|\lambda_i|)/2$ & $0.698$  \\
				$\max(|\lambda_i|)$ & $0.121$ \\
                $2(\max(|\lambda_i|))$ & $8.6\times 10^{-6}$ \\
				\hline
			\end{tabular}
	\end{center}
\end{table}

For our sixth experiment, we wish to obtain a real $13\times 13$ symmetric matrix with eigenvalues 1, 2, 2, 2, 3, 3, 3, 4, 4, 5, 5, 6 e 7, that is,
 having multiplicities 1, 3, 3, 2, 2, 1, and  1, and the following pattern:

\[ \left[ \begin{array}{ccccccccccccc}
	a & f & f & g  & h & 0 & 0 & 0 & 0 & 0 & 0 & 0 & 0 \\
	f & 4 & 0 & 0 & 0 & 1 & 1 & 0 & 0 & 0 & 0 & 0 & 0 \\
	f & 0 & 4 & 0 & 0 & 0 & 0 & 1 & 1 & 0 & 0 & 0 & 0  \\
	g  & 0 & 0 & b & 0 & 0 & 0 & 0 & 0 & i & i & 0 & 0\\
	h & 0 & 0 & 0 & c & 0 & 0 & 0 & 0 & 0 & 0 & j & k \\
	0 & 1 & 0 & 0 & 0 & 3 & 0 & 0 & 0 & 0 & 0 & 0 & 0\\
	0 & 1 & 0 & 0 & 0 & 0 & 3 & 0 & 0 & 0 & 0 & 0 & 0\\
	0 & 0 & 1 & 0 & 0 & 0 & 0 & 3 & 0 & 0 & 0 & 0 & 0\\
	0 & 0 & 1 & 0 & 0 & 0 & 0 & 0 & 3 & 0 & 0 & 0  & 0\\
	0 & 0 & 0 & i & 0 & 0 & 0 & 0 & 0 & 4 & 0 & 0  & 0\\
	0 & 0 & 0 & i & 0 & 0 & 0 & 0 & 0 & 0 & 4 & 0  & 0\\
	0 & 0 & 0 & 0 & j & 0 & 0 & 0 & 0 & 0 & 0 & d  & 0\\
	0 & 0 & 0 & 0 & k & 0 & 0 & 0 & 0 & 0 & 0 & 0  & e \\
\end{array} \right], \] \\ [2mm]
where  $a, b, c, d, e, f, g, h, i, j$, and $k$  are the considered variables, and they must belong to the nonzero set of indices  $I_{nz}$.
 This problem was posed as a challenging one in \cite[Section 7]{sutton}.
We will apply GLODS using different pairs of bounds for the variables, including the following ones:
$low_1=0.3\times ones(11,1)$,  $up_1=6\times ones(11,1)$; $low_2=0.5\times ones(11,1)$,
$up_2=5\times ones(11,1)$; $low_3=0.8\times ones(11,1)$,   $up_3=4.5\times ones(11,1)$; and
$$ low_4=( 0.5,  0.5,  0.5,  0.5,  0.5,  0.5, 0.5,  0.5, -5, -5, -5)\; \mbox{ and }
\; up_4=(5, 5, 5, 5, 5, 5, 5, 5, -0.5, -0.5, -0.5). $$

Solving our optimization model with  GLODS  produces a variety of solutions, all of them meeting the structural
 requirements. We now present four of the obtained matrices:   \\ [2mm]
  Exp6-Matrix-1 is completed with \\
 {\footnotesize  $(a,b,c,d,e,f,g,h,i,j,k)=(4.20377, 4.05002, 4.15003, 2.38052, 4.21566, 1.40729, 0.5501, 1.60692, 1.43179, 0.86127, 0.66676)$}, \\  [2mm]
  Exp6-Matrix-2  with \\
  {\footnotesize  $(a,b,c,d,e,f,g,h,i,j,k)=(4.3053, 3.0001, 4.2924, 4.7185, 2.6838, 0.7028,  1.7404, 0.9992,
  	1.0000, 2.0689, 0.7006)$}. \\  [2mm]
   Exp6-Matrix-3  with \\
  {\footnotesize  $(a,b,c,d,e,f,g,h,i,j,k)=(4.2608, 4.3027, 3.4959, 2.5937, 4.3468, 1.1568,  0.9950, 1.5778, 1.5175, 0.8287, 0.8922)$}. \\  [2mm]
    Exp6-Matrix-4  with \\
  {\footnotesize  $(a,b,c,d,e,f,g,h,i,j,k)=(4.5072, 4.3776, 3.1301, 4.0992, 2.8858, 0.8439,  0.9746, 1.8109, -1.542, -1.2384, -0.595)$}. \\  [2mm]
  We note that the matrices are different, but they all have the  given  eigenvalues
  {with the required precision}. In  Table \ref{Exp6}  we report the behavior
of GLODS to produce the 4 different solutions described above.

\begin{table}[H]
	\caption{Performance of GLODS for the $13\times 13$  symmetric matrix and  the considered bounds.}  \label{Exp6}
 	\begin{center}
		\begin{tabular}{cccccc} \hline
	Solution &		Bounds  &  Iter &  sum$(|NZ|)$  &  NumEvalf & $\min(|NZ|)$  \\ \hline
	Exp6-Matrix-1 &		$low_1\; \leq a,\ldots, k\leq \;up_1$  & 291 &  $34.86$ & 4182 &  0.5501 \\
	Exp6-Matrix-2 &		$low_2\; \leq a,\ldots, k\leq \;up_2$ & 671 &  $36.83$ & 9080 &  0.7006 \\
	Exp6-Matrix-3 &		$low_3\; \leq a,\ldots, k\leq \;up_3$ & 416 &  $38.27$ & 6230 &  0.8287   \\
	Exp6-Matrix-4 &		$low_4\; \leq a,\ldots, k\leq \;up_4$ & 430 &  $37.78$ & 6749 &  0.5950  \\
   \hline
		\end{tabular}
	\end{center}
\end{table}

	Now, fixing the   bounds as follows: $low_2\; \leq a,\ldots, k\leq \;up_2$, in
	Table \ref{Exp6tau} we report the value of
	$\|vec(\Lambda(X) - \Lambda)\|_2$  for an increasing set of values of the parameter
	$\tau>0$. It can be observed that for small values of $\tau$ the obtained matrix
	fails to reproduce the given eigenvalues {with the required precision}, and also
	 that as  $\tau$ increases the obtained
	eigenvalues get closer to the given ones. As reported in Table  \ref{Exp6}, for $\tau = 2(\max(|\lambda_i|))$,
	the obtained matrices reproduce the given eigenvalues {with the required tolerance}.

\begin{table}[H]
	\caption{Performance of GLODS for  $13\times 13$  symmetric matrix,
		using the bounds $low_2 \leq I_{nz} \leq up_2$, and  several values of the parameter $\tau$.}  \label{Exp6tau}
	 	\begin{center}
			\begin{tabular}{cc} \hline
				$\tau$  &  $\|vec(\Lambda(X) - \Lambda)\|_2$  \\ \hline
				$\max(|\lambda_i|)/10$ &  $ 5.06$   \\
				$\max(|\lambda_i|)/5$ &  $1.28$ \\
				$\max(|\lambda_i|)/2$ & $0.165$  \\
				$\max(|\lambda_i|)$ & $5.2\times 10^{-5}$ \\
                $2(\max(|\lambda_i|))$ & $9.1\times 10^{-6}$ \\
				\hline
			\end{tabular}
	\end{center}
\end{table}

  In summary, we observe that GLODS successfully solves the considered problems by setting
 $\tau>0$  to twice the largest absolute value of the given eigenvalues. This yields matrices that satisfy the
  required structural constraints in all cases. We also observe that it was consistently convenient to maintain both the simple bounds that define $\Omega$ and the second term of the function $f(X)$, i.e.,
  $- \sum_{(i,j)\in I_{nz}}  \ln(|x_{ij}|)$. This is because keeping both yields solutions for which the
   non-zero inputs have a higher average absolute value, as well as a higher minimum non-zero value.

\subsection{Performance of heuristic strategies}

We now report the performance of the heuristic strategies, described in Section \ref{Heurist}, on the same
 problems for which we reported  the performance of GLODS  in Subsection \ref{glods}.   For the first
 two heuristic strategies considered (GA and PSO),  we are using the MATLAB built-in implementations available
  in its  Optimization Toolbox. For the CMA-ES strategy we use the implementation in MATLAB freely available at
 {\bf https://github.com/CMA-ES/CMA-ES.github.io/blob/master/cmaes.m}. \\

For all experiments, as in Subsection \ref{glods}, we set  $\tau=2(\max(|\lambda_i|))$.
  In all cases, we considered two different sets of lower and upper bounds on the variables in $I_{nz}$: $[0.5,5]$ and $[0.5,10]$.
   We also considered negative intervals for the set $I_{nz}$, but the heuristics consistently performed less effectively than when
   dealing with positive intervals, and hence, we are not reporting them.
    The variables in  $I_{free}$ were bounded in $[-M, M]$ for $M>0$. We explored the behavior of the three heuristics for different values
    of $M$, from $10^{1}$ to $10^{5}$, and the results were, in most cases, quite similar or slightly worse. However, in a few cases,
    reducing the intervals increased the number of evaluations considerably.  Therefore, we have decided to set
     $M=10^{5}$.
    We note that regardless of the number of iterations or function  evaluations allowed, there were some problems
     for which some of the heuristics were not able to  reproduce the given eigenvalues {with the
  	 required tolerance}, only reaching a decimal precision. {Runs that did not reach the target eigenvalue tolerance (defined as $\|vec(\Lambda(X) - \Lambda)\|_2\leq 10^{-5}$)
    are not included in the following tables.}
 Additionally, as  in  Subsection \ref{glods}, we also examined the behavior of the heuristics when the second term of
  the objective function is omitted, i.e., when $f(X) = \|vec(\Lambda(X) - \Lambda)\|_2$, and the results were
  quite similar, except a higher failure rate was observed in several cases. Due to space limitations, we only
  report the results when both terms of the  function $f$ are considered.  \\

 Due to the random nature of heuristic algorithms, for each experiment and each set of bounds, we perform ten 
 runs {(with ten different seeds)}  and report the one corresponding to the smallest number of
  function  evaluations among those that successfully obtain the eigenvalues {with the
  	 required precision}.   In each case, we report the sum {in absolute value}  of all the
  entries in $I_{nz}$ that appear in the solution matrix (sum$(|NZ|)$),
  and also the smallest entry {in absolute value}  in $I_{nz}$  ($\min(|NZ|)$). We also report
the required number of function evaluations (NumEvalf) and the average number of function evaluations (AvEvalf)
 of the successful runs, i.e., those in which the given eigenvalues are obtained. When the number of successful runs is strictly less than 10,
 we put that number in parenthesis to the right of NumEvalf.  \\

In Table \ref{Heuristics1} we present the performance of the Genetic Algorithm  to find the symmetric matrices
reported in all experiments. To ensure a good exploration of the search space we considered a population 
size of 50 and a maximum  of 1000 iterations. {The initial population is randomly generated 
	 within the bounds with a uniform distribution.}

 \begin{table}[H]
  	\caption{Performance of the Genetic Algorithm for the symmetric matrices experiments.	}  \label{Heuristics1}
 \begin{center}
	\small
	\begin{tabular}{ccrrrrr} \hline
		Experiment            & $n$ &Bounds   &   NumEvalf & AvEvalf &sum$(|NZ|)$& $\min(|NZ|)$  \\ \hline
		\multirow{ 2}{*}{Exp1}&	& [0.5,5]                            &   6398  & 9260 & 5.98  &   0.9387\\
		                                   &  & [0.5,10]                          &    7150 & 9956 & 6.07  & 0.9448 \\ \hdashline
		\multirow{ 2}{*}{Exp2}&	& [0.5,5]                             &  16268   & 36680    & 11.49	&   0.5355   \\
		                                   &  & [0.5,10]                          &  16691   &  29569  & 9.65  &  0.5724     \\ \hdashline
		\multirow{ 2}{*}{Exp3}&	& [0.5,5]                             &  4941   &  8114  & 14.12  &  1.2975  \\
		                                   &                                            & [0.5,10]    & 5552    & 7470   & 14.09 & 1.2641    \\ \hdashline
		 \multirow{ 6}{*}{Exp4}&\multirow{ 2}{*}{4}	& [0.5,5] &  6116   & 9176   & 5.94  & 0.9560   \\
		                                       &                           & [0.5,10]&  6868   &  9312  & 5.90  &  0.9377   \\ \cdashline{2-7}
		                                       &\multirow{ 2}{*}{7}& [0.5,5] & 14764    &  19318  & 9.49  & 0.6250   \\
		                                       &                           & [0.5,10] & 14717    &  19784  & 9.94 &  0.5226	   \\ \cdashline{2-7}
		                                       &\multirow{ 2}{*}{10}& [0.5,5] &  28770{\tiny (9)}   & 36467   & 16.50  & 0.5626   \\
		                                       &                               & [0.5,10] &   29146  &  35026  & 17.44 &  0.7047   \\ \hdashline
		  \multirow{ 2}{*}{Exp5}&	                        & [0.5,5]    &  4753   &  11850  & 9.97  &  0.5457  \\
		                                       &                               & [0.5,10] &  5881   &  10731  & 9.84 &  0.5165   \\ \hdashline
		  \multirow{ 2}{*}{Exp6}  &	                         & [0.5,5]  &  23506{\tiny (1)}   &  23506  & 37.68  &  0.5695  \\
		                                       &                                & [0.5,10]    & 25621{\tiny (2)}    &  35585  & 35.89 &   0.5467  \\
		\hline
	\end{tabular}
\end{center}
  \end{table}

In Table \ref{Heuristics2} we report the performance of the Particle Swarm heuristic on all the experiments.
 To ensure a good exploration of the search space we considered 50 particles and a maximum of 1000 iterations.
{The positions of the particles are randomly generated with a uniform distribution 
	within their bounds.}

\begin{table}[H]
\caption{Performance of the Particle Swarm Optimization scheme for the symmetric matrices of all the experiments.}  \label{Heuristics2}
  \begin{center}
	\small
	\begin{tabular}{ccrrrrr} \hline
		Experiment            & $n$ &Bounds   &   NumEvalf & AvEvalf &sum$(|NZ|)$& $\min(|NZ|)$  \\ \hline
		\multirow{ 2}{*}{Exp1}&	                             & [0.5,5]      & 9200    & 14735  & 5.68  & 0.8177  \\
		                                   &          & [0.5,10]    & 8500{\tiny (8)}  & 14781 &  6.07 & 0.7448 \\ \hdashline
		\multirow{ 2}{*}{Exp2}&	                  & [0.5,5]      & 15050{\tiny (9)}    & 25761  & 11.79 & 0.8820   \\
		                               &       & [0.5,10]    &  16800{\tiny (8)}   &  22956  & 10.80  &  0.5293   \\ \hdashline
		\multirow{ 2}{*}{Exp3}&	       & [0.5,5]      &  4600{\tiny (8)}   & 10256   &  14.13 & 1.3305   \\
		                                &           & [0.5,10]    &  5800{\tiny (6)}   & 7133   & 14.14 &  1.3498   \\ \hdashline
		 \multirow{ 6}{*}{Exp4}&\multirow{ 2}{*}{4}	& [0.5,5]      &  9000   &  14360  &5.94   &  0.8990  \\
		                                    &                             & [0.5,10]    &   9150{\tiny (8)}  &  13044  & 5.51 &  0.7255   \\ \cdashline{2-7}
		                                    &\multirow{ 2}{*}{7}	& [0.5,5]      &  15100{\tiny (9)}   &  23406  &  11.51 & 0.7947   \\
		                                    &                              & [0.5,10]    &  15400{\tiny (9)}   & 20189   & 10.04 &  0.5007   \\  \cdashline{2-7}
		                                     &\multirow{ 2}{*}{10}& [0.5,5]     & 20700{\tiny (9)}    &  27250 & 16.38  &  0.6232  \\
		                                    &                               & [0.5,10]    & 24400    &  31645  & 17.09 & 0.6987    \\ \hdashline
		  \multirow{ 2}{*}{Exp5}&	                     & [0.5,5]      &  2750{\tiny (8)}   & 7163   &  9.78 & 0.5000   \\
		                                   &                                & [0.5,10]    &  2900{\tiny (9)}   & 6344   &  9.78 & 0.5000   \\ \hdashline
		  \multirow{ 2}{*}{Exp6}&	       & [0.5,5]      &  11350{\tiny (5)}   &   36600 &  37.83 &  0.6490  \\
		                                   &     & [0.5,10]    &  11000{\tiny (6)}   &  19058  & 37.20 &  0.5866   \\
		\hline
	\end{tabular}
\end{center}
 \end{table}

Table \ref{Heuristics4} presents the results of applying the CMA-ES heuristic to all the experiments, considering an initial point uniformly
 distributed within the corresponding bounds.

\begin{table}[H]
\caption{Performance of CMA-ES for the symmetric matrices of all the experiments, except for Experiment 4, $n = 10$, and
the bounds for  $I_{nz}$ in $[0.5, 10]$, since the given eigenvalues could not be obtained in any of the 10 runs. }  \label{Heuristics4}
 \begin{center}
 	\small
	\begin{tabular}{ccrrrrr} \hline
		Experiment            & $n$ &Bounds   &   NumEvalf & AvEvalf &sum$(|NZ|)$& $\min(|NZ|)$  \\ \hline
		\multirow{ 2}{*}{Exp1}&	                             & [0.5,5]      &  4609   & 35679 & 5.81  & 0.5071  \\
		                                   &                & [0.5,10]    &  12169  & 46523 & 4.10  &  0.6135\\ \hdashline
		\multirow{ 2}{*}{Exp2}&	              & [0.5,5]      &  12420{\tiny (8)}  &   137020  & 10.31 & 0.8664   \\
		                                   &    & [0.5,10]    &  7943   & 117924   & 7.60  &  0.6499   \\ \hdashline
		\multirow{ 2}{*}{Exp3}&	               & [0.5,5]      &  5017   &  7360  & 14.14  &  1.4142  \\
		                                   &    & [0.5,10]    &  6673   &  7879  & 14.14 &  1.4142   \\ \hdashline
		 \multirow{ 6}{*}{Exp4}&\multirow{ 2}{*}{4}	& [0.5,5]      &  22033   & 39035   & 5.99  & 0.9368   \\
		                                    &                             & [0.5,10]    &  3241   &  35960  & 5.18 &  0.6453   \\ \cdashline{2-7}
		                                    &\multirow{ 2}{*}{7}	& [0.5,5]      &  104292{\tiny (5)}    &  215458  & 11.89  & 0.8092   \\
		                                    &                              & [0.5,10]    &   62789{\tiny (5)}  & 234266   & 11.98 &   0.9104  \\ \cdashline{2-7}
		                                     &                       10   & [0.5,5]     & 525949{\tiny (3)}    &  658721  & 17.98  & 0.9763   \\ \hdashline
		  \multirow{ 2}{*}{Exp5}&	                     & [0.5,5]      &  1555   &  2303  & 9.78  &  0.5000  \\
		                                   &                                & [0.5,10]    &  2164   & 2583   & 9.78 &  0.5000   \\ \hdashline
		  \multirow{ 2}{*}{Exp6}&	           & [0.5,5]      &   65572{\tiny (7)} & 99682   &  36.86 &   0.7579 \\
		                                   &  & [0.5,10]    &   51921{\tiny (8)}  &  160914  & 37.24  & 0.5000    \\
		\hline
	\end{tabular}
 \end{center}
 \end{table}

Summing up, we observe that  GA, PSO and CMA-ES successfully solve  the  considered problems obtaining matrices that satisfy the
 required constraints.
We note that, in some cases, CMA-ES requires a significant increase in the number of function evaluations.
However, for cases such as  Exp4 and Exp5 ($n=4$ and $I_{nz} \in [0.5, 10]$), CMA-ES requires the fewest function evaluations to obtain
values of   $\min(|NZ|)$ and  sum$(|NZ|)$ that are quite similar to the ones obtained by GA and PSO.

\section{Final comments and perspectives} \label{conclu}

{This paper introduces a novel, variable-efficient derivative-free optimization (DFO) modeling
	 framework for structured inverse symmetric matrix problems with prescribed eigenvalues. Our approach
	 minimizes a Lipschitz continuous objective function over a convex set, significantly reducing the number
	  of variables compared to manifold optimization methods (e.g., only 5 variables in our third experiment
	   vs. 100 in \cite{chehab}). We present the first comprehensive empirical comparison of deterministic
	   (GLODS) and heuristic DFO solvers for  this model. Both families produce high-quality solutions for a wide
   variety of matrices of different sizes and patterns. GLODS offers a superior balance of efficiency and
    solution quality, especially in terms of the absolute value of the smallest non-zero entry ($\min(|NZ|)$)
     and the sum of the absolute values of all non-zero entries (sum$(|NZ|)$). While heuristics are viable,
      they require substantially more function evaluations than GLODS. Our novel DFO approach outperforms
     \cite{chehab} in terms of both variable count and solution quality, demonstrating its  effectiveness and
      efficiency for this class of problems. \\

 Finally, as we look ahead, several direct perspectives emerge for extending this work. These include
   scaling and adapting the proposed approach to larger, sparse matrices (where $n >100$), exploring the
   inclusion of more complex constraints for the $I_{nz}$ set, and investigating the possible
    integration of  more advanced single- or multi-objective DFO algorithms, particularly those tailored to
     non-smooth problems, such as model-based methods within this framework; see, e.g.,
     \cite{AudetHare,AudetHare20, bras20, ConnEtAl, Custodio24, Custodio18, Dzahini, Garman, Larson}.
      Each of these directions offers  promising ways to further extend the applicability of the proposed approach.}  \\ [2mm]

\noindent

\newpage\noindent
 {\bf Funding.} \\
    This work is funded by national funds through the FCT -   Funda\c c\~ao para a Ci\^encia e a Tecnologia, I.P.,
    under the scope of projects UID/00297/2025  (https://doi.org/10.54499/UID/00297/2025) and UID/PRR/00297/2025
    (https://doi.org/10.54499/UID/PRR/00297/2025) (Center for Mathematics and Applications - NOVA Math). 
    The second author received financial support from CNPq under contract no. 409837/2025-3. \\ [2mm]

   \noindent
  {\bf Conflict of interest.} No potential conflict of interest was reported by the authors. \\ [2mm]

  \noindent
   {\bf Data availability.}  The codes and datasets generated during and/or analyzed during the
 	 current study are available at https://github.com/EvelinKrulikovski/DFOapproach. \\ [2mm]

 \noindent
{\bf Consent for publication.} All authors read and approved the submission of the final manuscript.  \\ [2mm]

\noindent
{\bf Author's contributions.} The three authors contributed equally to the  study, conception, design, data collection,
conceptualization, methodology, and analysis.

\end{document}